\documentclass[11pt]{amsart}
\usepackage{amscd,amsmath,amssymb}
\usepackage[mathscr]{eucal}
\usepackage{graphicx}
\newtheorem{thm}{Theorem}[section]

\newtheorem{cor}[thm]{Corollary}
\newtheorem{lem}[thm]{Lemma}
\newtheorem{prop}[thm]{Proposition}
\newtheorem{defn}[thm]{Definition}
\newtheorem{rem}[thm]{Remark}

\numberwithin{equation}{section}

\newcommand{\abs}[1]{\left\vert#1\right\vert}

\newcommand{\eps}{\varepsilon}

\begin{document}
\title{Parabolic Foliations on 3-manifolds.}
\author{Vladimir Krouglov.}
\maketitle

\section{Introduction.}
It is well known that every closed orientable three manifold admits a foliation. This statement becomes false if we consider additional restrictions on the geometry and topology of the leaves. For example, it is well known that foliations by the minimal surfaces do not exist on a three-sphere (relative to any metric). Analogously, the classes of totally umbilical foliations and totally geodesic foliations also exist not on every three-manifold.

A. Borisenko introduced new classes of foliations on Riemannian manifolds which have restrictions on the extrinsic geometry of the leaves, namely elliptic, parabolic and strong saddle (or hyperbolic) foliations. The study of existence of these foliations on 3-manifolds was initiated by D.Bolotov in \cite{Bol}. In this work, among the other results, he defines a metric on a solid torus such that Reeb component is a parabolic foliation. In \cite{Bol2} he gives examples of strong saddle foliations on torus bundles over the circle and on a three sphere. In particular, a foliation in a Reeb component is not a topological restriction to the existence of strong saddle foliations. In \cite{Kru1}, author showed that in fact every closed orientable three manifold admits a strong saddle foliation.

It is well known that closed orientable 3-manifolds do not admit elliptic foliations. Namely, their existence contradicts to an integral formula $\int_{M} H = 0$ (here  $H$ -- stands for the mean curvature) since if it is elliptic with respect to some metric its total mean curvature cannot be zero.

The last open problem was the existence of parabolic foliations on closed orientable 3-manifolds. In this paper we give positive answer on this question.
\begin{thm}
Every closed orientable 3-manifold admits a parabolic foliation.
\end{thm}

Note, that there are no parabolic foliations on $S^3$ with respect to a standard metric.

Parabolic foliations of codimension greater than one were studied in \cite{Bor}.

This paper is organized as follows: in Section 2 we recall some definitions and constructions from the topology of foliations on 3-manifolds. In Section 3 we construct several local models of parabolic foliations. In Section 4 we define a parabolic foliation on a 3-spehere which is a turbularization of a Reeb foliation along an arbitrary knot. In Section 5 we show how to perform a Dehn surgery on this knot to obtain a parabolic foliation on every closed orientable $3$-manfiold.

\textbf{Acknowledgements:} I would like to express my gratitude to Dmitry Bolotov for his constant support and attention to this work. I would also like to thank Professor Alexander Borisenko for his help and valuable advices during the preparation of the paper.
\section{Basic Definition.}
\subsection{Foliations on three-manifolds.}

In this section we will recall basic necessary about the foliations on three-manfiolds.

Let $\mathcal{F}$  be a foliation on a closed three-manifold. It defines a two-dimensional distribution of planes tangent to the leaves. However, not every plane distribution defines a foliation. A distribution is called integrable if it defines a foliation. Classical Frobenius theorem gives necessary and sufficient conditions for a distribution to be integrable. Here we will recall a three-dimensional version of this theorem.
\begin{thm}(Frobenius)
A distribution of planes $\xi$ on a three-manifold is integrable if and only if for every pair of local sections $X$ and $Y$ of $\xi$ its Lie bracket belongs to $\xi$.
\end{thm}
Recall, that a distribution is called transversally orientable if there is a globally defined vector field transverse to it. In this case there is a globally defined one-form  $\alpha$ such that $Ker(\alpha)_p = T_p L$ (where $L$ is a leaf through $p$). It is easy to rewrite conditions of Frobenius theorem in terms of the form $\alpha$: a distribution is integrable if and only if $\alpha \wedge d\alpha = 0$. \\ \\
\textbf{Example 2.2:}\emph{Reeb foliation on $D^2 \times S^1$}. \\
Consider the following $C^\infty$-smooth function on  $[0,1]$:
\begin{enumerate}
\item{The function $f$ is a smooth increasing function on $[0,1]$.}
\item{There is an $\eps>0$ such that for any $x\in [0, \eps)$ the value of $f(x)$ is equal to zero and for any $x \in (1-\eps, 1]$, $f(x)=1$.}
\end{enumerate}
On the solid torus $D^2 \times S^1$ with the cylindrical coordinates $((r, \phi), t)$ define the following one-form:
$$
\alpha = f(r) dr + (1 - f(r)) dt
$$
From Frobenius theorem, a distribution of planes defined by the kernel of $\alpha$ is integrable since:
$$
\alpha \wedge d\alpha = (f(r) dr + (1 - f(r)) dt) \wedge (- f'(r) dr \wedge dt ) = 0
$$
Therefore $\alpha$ defines a foliation on $D^2 \times S^1$. We will denote this foliation by $\mathcal{F}_R$ and call it the Reeb foliation in a solid torus.
\begin{rem}
In the literature a Reeb foliation is usually defined as a foliation of $D^2 \times S^1 = \{((r, \phi), t): r \in [0, 1], \phi, t \in [0, 2\pi)\}$ by the levels of function $h(r, \phi, t) = (r^2 - 1) e^t$. It is obvious that foliation defined above is isotopic to this foliation. This justifies the name in Example 2.2.
\end{rem}
\subsection{Extrinsic geometry of foliations.}

Assume now that $M$ is a Riemannian manifold with a scalar product $g$ and associated Levi-Civita connection $\nabla$. Consider a foliation $\mathcal{F}$ on $M$. For each pair of vector fields $X$ and $Y$ on $M$ tangent to $\mathcal{F}$, define a second fundamental form of $\mathcal{F}$ with respect to a unit normal $n$ by
$$
B(X, Y) = g(\nabla_{X} Y, n)
$$
Using the scalar product in the tangent bundle we may define the following linear operator $A_n$.
$$
 B(X, Y) = g(A_n X, Y)
$$
This operator is called a Weingarten operator. Since $A_n$ is symmetric it has two real eigenvalues. These eigenvalues are the principal curvature functions. A product $K_e = k_1 k_2$ is called an extrinsic curvature of $\mathcal{F}$.

Depending on the sign of extrinsic curvature one can consider the following natural classes of foliations on three-manifolds.
\begin{defn}(Borisenko) A codimension one foliation on a three-manifold is called:
\begin{enumerate}
\item{Parabolic, if there is metric such that $K_e = 0$.}
\item{(Strong)saddle, if there is a metric such that $(K_e < 0)K_e \le 0$. }
\item{Elliptic, if there is a metric such that $K_e > 0$.}
 \end{enumerate}
\end{defn}
\begin{rem}
As already mentioned in the introduction, there are no elliptic foliations on closed oriented three-manifolds. Since $K_e > 0$, the functions of principal curvatures are nowhere zero and has to be simultaneously greater than zero or less than zero. From the integral formula
$$
0 = \int_M H = \frac{1}{2}\int_M (k_1 + k_2) \ne 0
$$
\end{rem}
\begin{rem}
Note, that many geometric classes of foliations fall into one of the introduced classes. It is easy to see that minimal foliations are saddle, totally umbilical foliations have  $K_e \ge 0$, and totally geodesic foliations are parabolic.
\end{rem}
\subsection{Knots and braids.}

Recall that a knot in $S^3$ or $\mathbb{R}^3$ is an image of a circle $S^1$ under some $C^\infty$-smooth regular embedding. Two knots $K_0$ and $K_1$
are called isotopic if there is a smooth family of embeddings $K(t) : S^1 \to S^3 (\mathbb{R}^3)$ such that $K(0) = K_0$ and $K(1) = K_1$.

Consider two sets of points $A = \{(i, 0, 0), i = 1, \ldots, n\}$ and $B=\{(i, 0, 1), i = 1,\ldots,n\}$ in
$\mathbb{R}^3$. A smooth embedded curve $\gamma(t)$ is called descending if its $z$-coordinate is a strictly decreasing function of parameter  $t$.

A topological braid $K$ with $n$ strings is a collection of $n$ disjoint descending curves in $\mathbb{R}^3$, which connect the points from the set $B$ with the points of $A$. We say that two braids are isotopic if there is a smooth family of braids connecting them.

Consider a group with the generators $\sigma_1, \sigma_2, \ldots, \sigma_{n-1}$ and relations $\sigma_i \sigma_{i + 1}\sigma_i = \sigma_{i+1} \sigma_i \sigma_{i +1}$ for all $i$ and also $\sigma_i \sigma_j = \sigma_j \sigma_i$, in the case when $|i - j| \ge 2$. This group is denoted by $B_n$, and called a group of algebraic braids. There is a one-to-one correspondence between the isotopy classes of topological braids and elements of $B_n$. We use convention that to a generator $\sigma_i$ corresponds a topological braid which frontal projection is represented at the left part of Figure 1.
Second possible intersection corresponds to the element $\sigma_i^{-1}$. Therefore, isotopy classes of topological braids may be represented as a products $K = \sigma_1^{\pm 1}\sigma_2^{\pm 1}\ldots\sigma_N^{\pm 1}$.

\begin{figure}
\includegraphics[width = 5cm]{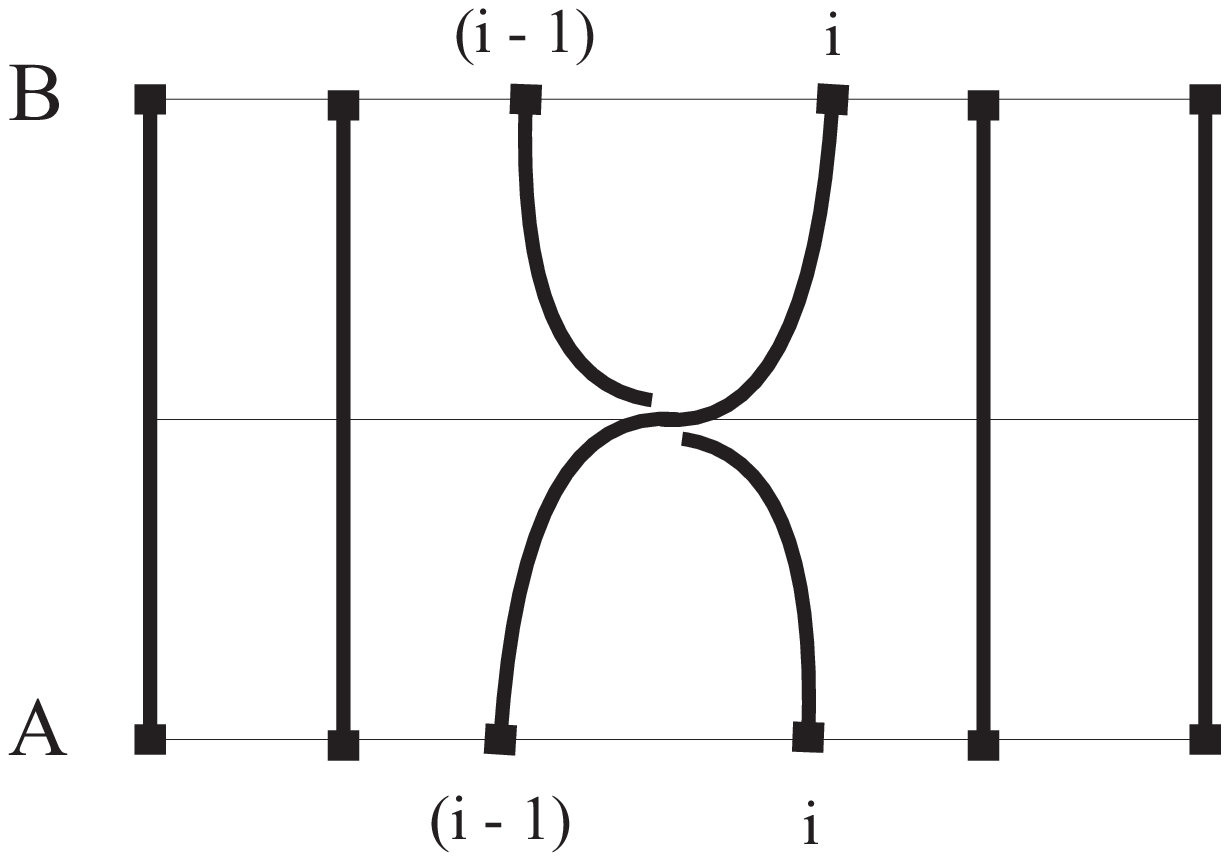} \ \ \ \ \
\includegraphics[width = 5cm]{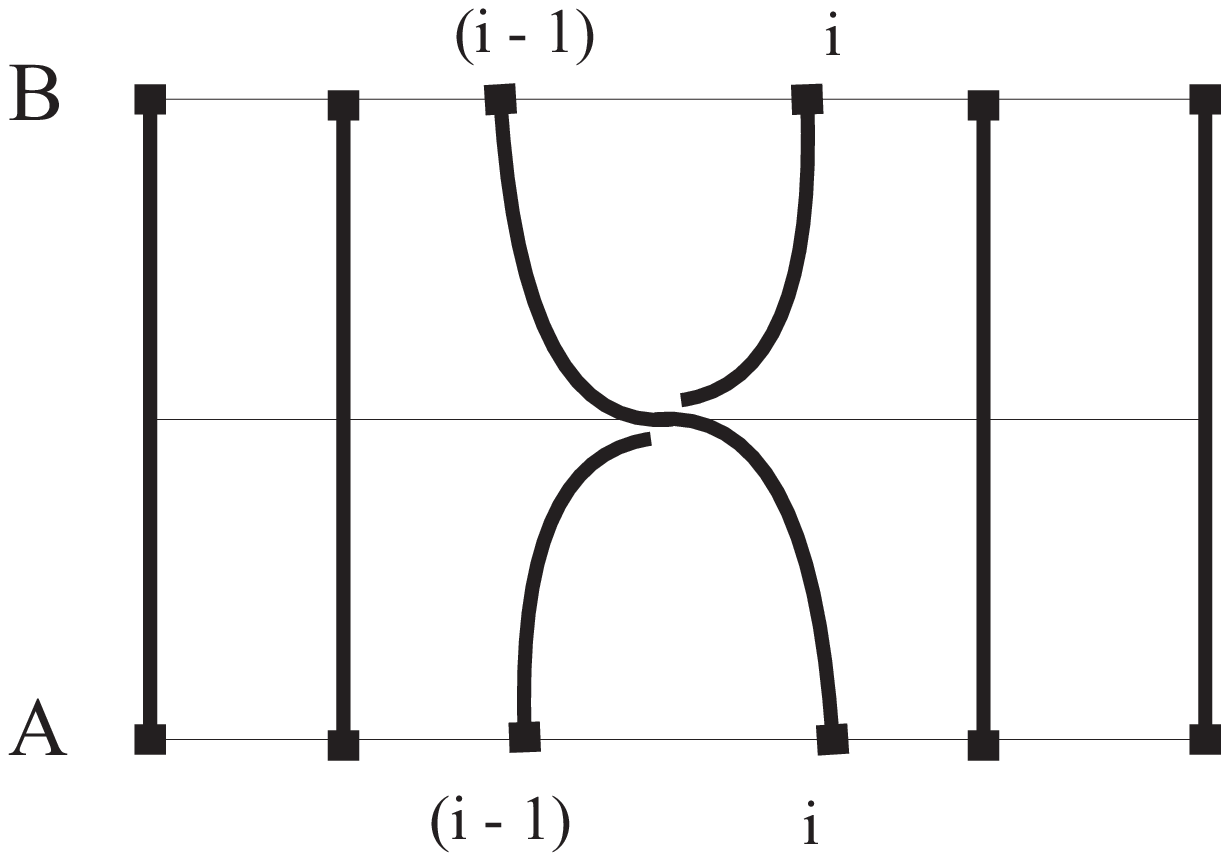}
\caption{Possible intersections in the frontal projection.}
\end{figure}

On the set of isotopy classes of topological braids define an operation of closure of a braid. A closure of a braid $K$ is a link (that is an embedded image of disjoint $S^1$) which is obtained from $K$ adding disjoint curves connecting $i$'th point of $A$ with $i$'th point of $B$.
 The following theorem holds(cf. \cite{PrSos}):
\begin{thm}
A map from the set of isotopy classes of topological braids which assigns to each braid its closure is surjective. In particular, each isotopy class of knots contains the closure of some braid.
\end{thm}
\subsection{Combinatorial presentation of three-manifolds.}
In this section we will give a sketch of proof that every closed orientable three-manifold admits a foliation.

Consider a knot $K$ in $S^3$. Let $N$ be some tubular neighborhood of $K$. Denote $X = \overline{S^3 \backslash N}$. Then $\partial X = \partial N = T^2$. Consider some homeomorphism
$$
h : \partial X \to \partial(D^2 \times S^1)
$$
and let $M = X \cup D^2 \times S^1 / (y \sim h(y), \  \mbox{for all $y \in \partial X$})$. It is easy to see that $M$ is a closed manifold.

This construction is called a Dehn surgery on a knot. Importance of this construction comes from the following theorem:

\begin{thm}\cite{PrSos} Every closed orientable manifold may be obtained by the Dehn surgery on some knot $S^3$.
\end{thm}

Recall the construction of a transversally orientable foliation on a closed orientable three-manifold. Consider a solid torus $D^2 \times S^1 = \{((r, \phi), t): r \in [0, 2], \phi, t \in [0, 2\pi)\}$ and let $\alpha = f(r) dr + (1 - f(r))dt$ where $f(r)$ is some smooth function on a segment $[0, 2]$, which satisfies the following conditions:
\begin{enumerate}
\item{$f(r)$ is a strictly increasing function on $[0, 1]$}
\item{$f(r)$ is a strictly decreasing function on $[1, 2]$}
\item{There is an $\epsilon$ such that for all $r \in (2-\epsilon, 2]$ the function $f(r) = 0$}
\item{$f(0) = 0$ and $f(1) = 1$}
\end{enumerate}
The form $\alpha$ defines some foliation on $D^2 \times S^1$. Denote it by $\mathcal{F}_T$.

It is obvious that $\mathcal{F}_T$ has a single compact leaf $\{r = 1\}$. $\mathcal{F}_T$ restricted on a solid torus $D^2(1) \times S^1 = \{(r, \phi, t) \in D^2 \times S^1: r\in[0, 1]\}$ is a Reeb foliation (see Example $2.2$).

It is well known that $S^3$ may be represented as a union of two solid tori, glued along the boundary torus. Gluing homeomorphism interchanges generators of the boundary torus. In each solid torus consider Reeb foliations $\mathcal{F}_R$. Since the gluing homeomorphism maps a leaf of the first Reeb component to the leaf of the second we see that the three-sphere admits a foliation which is the union of two Reeb components. We will also denote this foliation by $\mathcal{F}_R$.

Assume now that $K$ is a knot in $S^3$. From Theorem  $2. 4$ it is isotopic to the closure of some braid. We can further isotope this braid to make it everywhere transverse to the foliation of one of the solid torus by disks $D^2 \times \{t\}$. Since $\mathcal{F}_R$ is a foliation by disks in a small neighborhood of the core curve $r = 0$ we may assume that $K$ is transverse to $\mathcal{F}_R$. Cut a small tubular neighborhood of $K$ and glue back a solid torus with the foliation $\mathcal{F}_T$ inside. We will obtain a new foliation on $S^3$ which is a turbularization of the initial one along $K$.
Finally in order to obtain a foliation on $M$, cut a tubular neighborhood of $K$ up to the torus leaf and glue it back by the diffeomorphism of the boundary. It is easy to verify that since the boundary of this neighborhood is a leaf, the foliation is correctly defined on $M$. From Theorem $2.7$ every closed orientable three-manifold may be obtained this way, therefore every closed orientable three-manifold admits a foliation.
\subsection{Bump functions on $\mathbb{R}$.}

In the process of the proof we will often face with the following situation: in some finite segment $[a, b]$ we need to define a smooth function $f(t)$ in such a way that the following conditions are satisfied:
\begin{enumerate}
\item{$f(a) = f_0$}
\item{$f(b) = f_1$}
\item{There is an $\epsilon > 0$ such that:} \\
 - for all $t \in [a, a + \epsilon)$, $f(t) = f_0$ \\
 - for all $t \in [b - \epsilon, b)$, $f(t) = f_1$
\item{$f$ is monotone on $[a, b]$.}
\end{enumerate}
In the text of the paper we will refer to such functions as bump functions on $[a, b]$.
\begin{figure}
\includegraphics[width = 6.5cm]{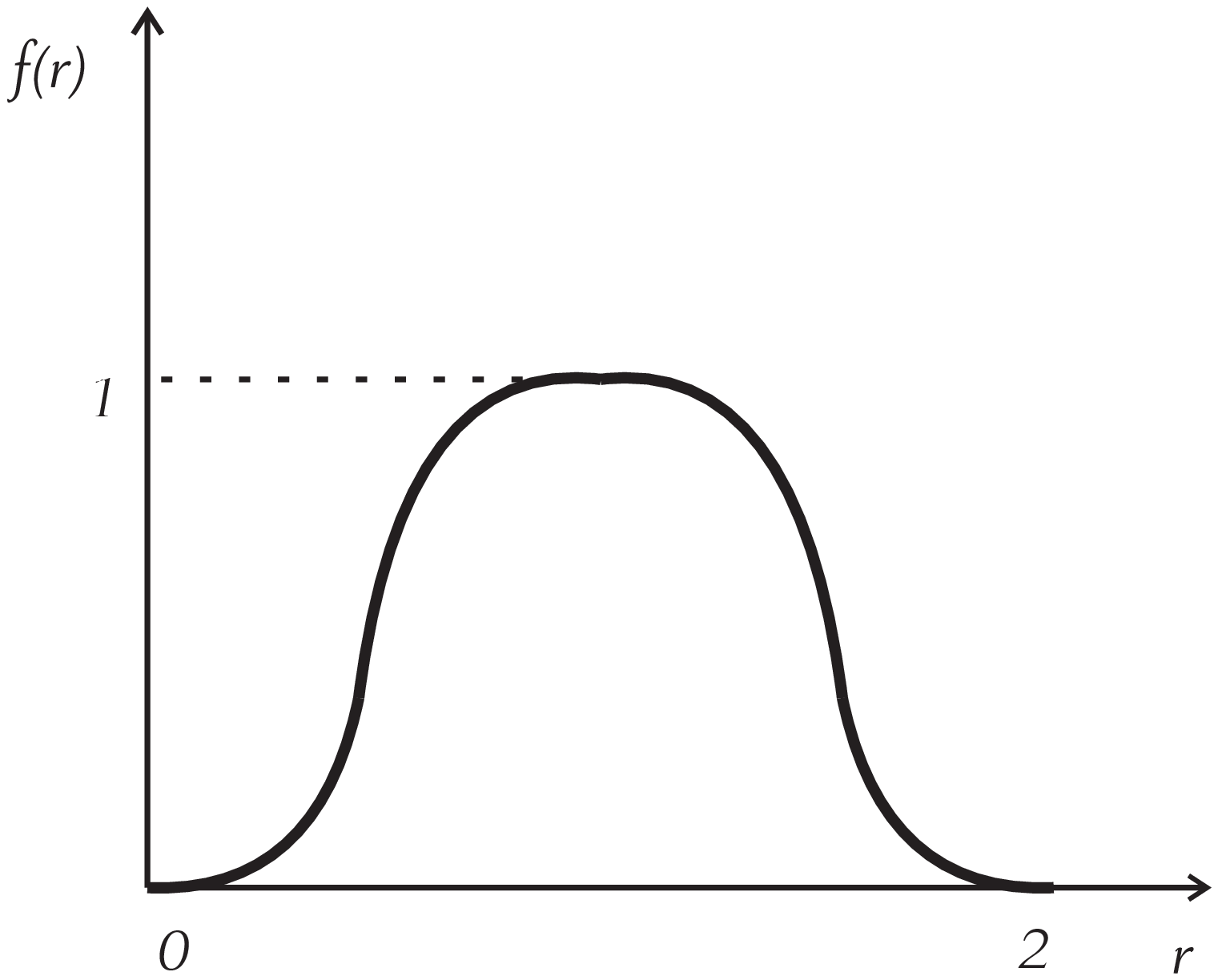}
\caption{Graph of $f(r)$ in the construction of turbularization.}
\includegraphics[width = 6.5cm]{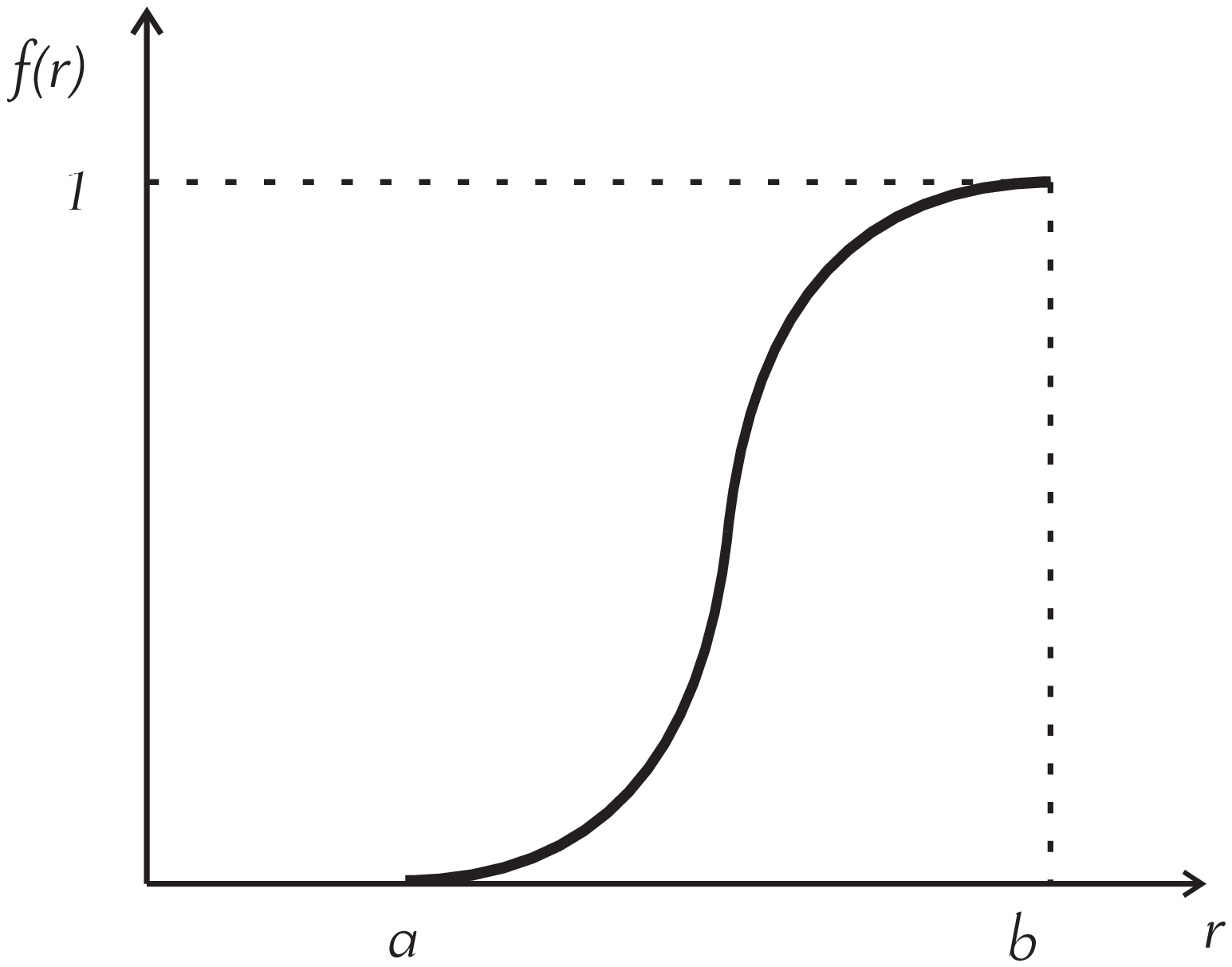}
\caption{Typical bump function on $[a, b]$.}
\end{figure}

For example, in the construction of Reeb component $\mathcal{F}_R$, $f$ is an increasing bump function on $[0, 1]$. Function $f(r)$ which arises in the construction of $\mathcal{F}_T$ is a union of two bump functions.
\section{Local models of parabolic foliations.}

In this section we will describe several local models of parabolic foliations on three-manifolds.
\subsection{Parabolic foliation on $\Sigma^2 \times [0,1]$.}

\begin{lem} Let $\Sigma^2$ be a compact parallelizable surface (possibly with the boundary).
Consider two Riemannian metrics $G$ and $H$ on $\Sigma^2$ that coincide in some neighborhood of the boundary $\partial \Sigma^2$.
Assume that $\mathcal{F}$ is a foliation of $M = \Sigma^2 \times [0, 1]$ by the surfaces $\Sigma^2 \times \{t\}$. Then, there is such Riemannian metric $g$ on $M$ that
\begin{enumerate}
\item{In some tubular neighborhood of $\Sigma^2 \times \{0\}$, $g = dt^2 + G(p)$, for all $p \in \Sigma^2$.}
\item{In some tubular neighborhood of $\Sigma^2 \times \{1\}$, $g = dt^2 + H(p)$, for all $p \in \Sigma^2$.}
\item{$\mathcal{F}$ is parabolic on $\Sigma^2 \times [0, 1]$ with respect to $g$}.
\item{There is a neighborhood $U$ of the boundary $\partial \Sigma^2$ such that for all $t \in [0,1], \   g(p, t)|_{U \times \{t\}} = G(p)$}
\end{enumerate}
\end{lem}
\emph{Proof:}Let  $\{ X_0, Y_0 \}$ be an orthonormal frame on $\Sigma^2$ with respect to $H$. The matrix of $G$ in this frame may be written as
$$
 G = G(p) = \left ( \begin{array}{ll} a(p) & b(p) \\ b(p) & c(p) \end{array} \right )
$$ for all $p \in \Sigma^2$.
Since $\Sigma^2$ is compact the functions $a$, $b$, and $c$ are bounded on $\Sigma^2$.

Consider a frame $( X, Y, n )$ on $M$ where $X = (X_0, 0)$, $Y = (Y_0, 0)$, and $n = \frac{\partial}{\partial t}$.
Let $N$ be such a neighborhood of the boundary $\partial \Sigma^2$ that $G|_N = H|_N$. Denote by $L = \overline{M \backslash N}$.

We are going to interpolate between $G$ and $H$ on $\Sigma^2 \times [0, 1]$ using for this the following Riemannian metric:
$$
g = g(p, t) = \left ( \begin{array}{ccc} a(p, t) & b(p, t) & 0 \\ b(p, t) & c(p, t) & 0 \\ 0 & 0 & 1  \end{array} \right )
$$
where $(p, t) \in \Sigma^2 \times [0, 1]$ and $a(p, t)$, $b(p, t)$, and $c(p,t)$ are function on $\Sigma^2 \times [0, 1]$. The matrix of $g$ is written with respect to a frame $(X, Y, n)$.

From the definition of $g$, $n$ is a unit normal vector field to $\mathcal{F}$. Calculate the matrix of the second fundamental form of the leaves relative to the normal $n$ in the basis $\{ X, Y \}$. From Koszul formula,
$$
2 g(\nabla_X X, n) = 2 X(g(X, n)) - n(g(X, X)) + g([X, X], n) - 2 g([X, n], X)
$$
where $\nabla$ is a Levi-Civita connection of $g$.
Note, that since $X$ is independent of $t$, $g([X, n], X) = 0$. Therefore $g(\nabla_X X, n) = - \frac{1}{2}\frac{\partial a}{\partial t}$.

Similarly,
$$
2 g(\nabla_Y Y, n) = 2 Y(g(Y, n)) - n(g(Y, Y)) + g([Y, Y], n) - 2 g([Y, n], Y) = - n g(Y, Y)
$$
Finally, we have:
$$
2 g(\nabla_X Y, n) =  X(g(Y, n)) + Y(g(X, n)) - n(g(X, Y)) + g([X, Y], n) -  g([Y, n], X)
$$
$$
- g([X, n], Y) = - n(g(X, Y)) + g([X, Y], n) -  g([Y, n], X) - g([X, n], Y)
$$
Since $\mathcal{F}$ is a foliation we have that $g([X, Y], n) = 0$. Using the fact that $X$ and $Y$ does not depend on $t$ we see that $g([Y, n], X)=0$ and $g([X, n], Y)=0$. Therefore $g(\nabla_X Y, n) = -\frac{1}{2} n(g(X, Y))$.
Consequently, we may conclude that the second fundamental form of the leaves is given by a matrix
$$
B =  - \frac{1}{2}\left ( \begin{array}{ll} \frac{\partial a}{\partial t} & \frac{\partial b}{\partial t} \\ \frac{\partial b}{\partial t} & \frac{\partial c}{\partial t} \end{array} \right )
$$
An extrinsic curvature of $\mathcal{F}$ with respect to $g$ equals to
 $$K_e = \frac{1}{4} \frac{ \frac{\partial a}{\partial t}\frac{\partial c}{\partial t} - \frac{\partial b}{\partial t}^2}{ac - b^2}$$
Take some subdivision $ 0 \le t_1 < t_2 < t_ 3 < t_4 < t_5 \le 1$ of the segment $[0, 1]$. Assume that $D$ is a positive real number which is greater than $\max_{p \in \Sigma^2}\{ a(p), c(p)\}$. We will choose an exact value of $D$ later in the proof.

Consider the following function $h$ on $\Sigma^2$
\begin{enumerate}
\item{$h(p) = 1$, for all $p \in L$}
\item{$h(p) = 0$ in some neighborhood of $\partial \Sigma^2$}
\item{$h$ is a smooth nonnegative function on $\Sigma^2$}
\end{enumerate}
Consider a function $\tilde{a}(p, t) = D f(t) + (1 - f(t)) a(p)$, where $f(t)$ -- is an increasing bump function on $[0,t_1]$ with $f(0) = 0 $ and $f(t_1) = 1$. Finally, let $$a(p, t) = h(p)\tilde{a}(p, t) + (1 - h(p))a(p)$$
On $[0, t_1]$ define the following matrix (with respect to a frame $(X, Y, n)$):
$$
g_D = g_D(p, t) = \left ( \begin{array}{ccc} a(p, t) & b(p) & 0 \\ b(p) & c(p) & 0 \\ 0 & 0 & 1  \end{array} \right )
$$
This matrix clearly defines a metric on $\Sigma^2 \times [0, t_1]$, since when $t = 0$ it has positive diagonal entries and $a(p, t)$ is nondecreasing on the segment $[0, t_1]$. From the definition of bump functions $g = dt^2 + G(p)$ in some tubular neighborhood of $\Sigma^2 \times \{0\}$. A foliation by the surfaces $\Sigma \times \{t\}$ is parabolic with respect to the introduced metric.

On the segment $[t_1, t_2]$ we may change $c(p)$ in the similar fashion. Consequently, on a segment $[0, t_2]$ we have:
\begin{enumerate}
\item{$g_D(p, 0) = \left ( \begin{array}{ll} G(p) & 0 \\ 0 & 1 \end{array} \right )$.}
\item{$g_D(p, t_2) = \left ( \begin{array}{lll} a(p, t_2) & b(p) & 0 \\ b(p) & c(p,t_2) & 0 \\ 0 & 0 & 1\end{array}\right )$.}
\item{$\mathcal{F}$ is a parabolic foliation on $\Sigma^2 \times [0, t_2]$ with respect to $g_D$.}
\end{enumerate}
Consider an increasing bump function $f(t)$ on a segment $[t_2, t_3]$ with $f(t_2) = 0$ and $f(t_3) = 1$. On $[t_2, t_3]$ define $g_D$ by the matrix:
$$
g_D = \left ( \begin{array}{ccc} a(p, t_2) + f(t)b(p) & b(p)(1 - f(t)) & 0 \\ b(p)(1 - f(t)) & c(p, t_2) + f(t)b(p) & 0 \\ 0 & 0 & 1  \end{array} \right )
$$
Since $b(p)$ is bounded on $\Sigma^2$ and is equal to zero on $N$, the matrix of $g_D$ is positively definite for some choice of $D$ (see the definition of $D$ earlier in the proof) and therefore defines a metric.

The foliation is parabolic with respect to the introduced metric since an extrinsic curvature is given by $K_e = \frac{1}{4 det(g_D)}(f'(t)^2 b(p)^2 - (-f'(t))^2 b(p)^2) = 0$.

On $\Sigma^2 \times \{ t_3 \}$ the matrix of $g_D$ is diagonal. Since we eliminated all non-diagonal elements of $g_D$ we may now freely decrease the diagonal elements of $g_D$. For this, consider decreasing bump function $k(t)$ on $[t_3, t_4]$ with $k(t_3) = 1$ and $k(t_4) = 0$. On $[t_3, t_4]$ define
$$
g_D = \left ( \begin{array}{ccc} k(t) a(p, t_3) + (1 - k(t)) & 0 & 0 \\ 0 & c(p, t_3) & 0 \\ 0 & 0 & 1 \end{array} \right )
$$
For all $t \in [t_3, t_4]$ the matrix $g_D$ is positively definite and therefore defines a metric. At $t = t_4$ the matrix of $g_D$ is given by:
$$
g_D = \left ( \begin{array}{ccc} 1 & 0 & 0 \\ 0 & c(p, t_3) & 0 \\ 0 & 0 & 1 \end{array} \right )
$$
Analogously, on $[t_4, t_5]$ we may decrease a diagonal element $c(p, t_3)$. Therefore we showed how to deform $g_D$ into a metric $\left ( \begin{array}{lll} 1 & 0 & 0 \\ 0 & 1 & 0 \\ 0 & 0 & 1 \end{array} \right )$, while keeping the parabolicity of $\mathcal{F}$. But this metric is nothing else as a matrix of $dt^2 + H(p)$ written with respect to a frame  $(X, Y, n )$.

Let $U$  be an open subset of $\{p \in \Sigma^2, \mbox{where}\  h(p) = 0\}$ which contains the boundary of $\Sigma^2$. On $[0, t_1]$ the function $a(p, t)|_U = a(p)|_U = 1$. The same holds for the function $c(p, t)|_U$ on a segment $[t_1, t_2]$. On $[t_2, t_3]$, since $b(p)|_U = 0$, the matrix of $g|_U$ is an identity matrix. On the segments $[t_3, t_4]$ and $[t_4, t_5]$, since $a(p, t)|_U = c(p, t)|_U = 1$ the matrix $g|_U$ is also an identity matrix. This finishes the proof of the lemma.

\begin{cor}
Consider the manifold $M = T^2 \times [0,1]$. Let $G(x, y)$ and $H(x, y)$ be any metrics on $T^2$. Then there is a metric $g$ on $M$ such that
\begin{enumerate}
\item{Foliation $\mathcal{F}$ by the tori $T^2 \times \{pt\}$ is parabolic.}
\item{The matrix of $g|_{T^2 \times \{0\}} = \left ( \begin{array}{cc} G(x, y) & 0 \\ 0  & 1\end{array}\right )$ and the matrix of $g|_{T^2 \times \{1\}}  =\left ( \begin{array}{cc} H(x, y) & 0 \\ 0  & 1\end{array}\right )$}
\end{enumerate}
\end{cor}
\subsection{Parabolic foliation on a solid torus.}
The following proposition is due to D.Bolotov.
\begin{lem}(D.Bolotov, \cite{Bol})
There is a foliation $\mathcal{F}$ and a metric $g$ on $D^2 \times S^1$ such that:
\begin{enumerate}
\item{$\mathcal{F}$ is parabolic with respect to $g$.}
\item{The foliation $\mathcal{F}|_{D^2(\frac{1}{3}) \times S^1}$ is a foliation by the totally geodesic disks $D^2(\frac{1}{3}) \times \{t\}$ and the foliation $\mathcal{F}|_{([\frac{2}{3}, 1] \times S^1) \times S^1}$ is a foliation by the totally geodesic tori $\{ r \} \times S^1 \times S^1$.}
\end{enumerate}
\end{lem}
\emph{Proof:} Consider the solid torus $D^2 \times S^1$ with the following coordinates on it
$$
D^2 \times S^1 = \{ ((r, \phi), t) : r \in [0, 1], \phi, t \in [0, 2\pi) \}
$$
Define  the one-from $\alpha$ on $D^2 \times S^1$ as:
$$
 \alpha = f(r)dr + (1 - f(r)) dt
$$
where $f(r)$ is such a smooth function on $[0,1]$ that
$$
 f(r) = \left \{ \begin{array}{l} 0, \ r \in [0, \frac{1}{3}] \\
                                  \mbox{is a strictly increasing function when} \   r \in (\frac{1}{3}, \frac{2}{3}] \\
                                  1, \   r \in (\frac{2}{3}, 1]
  \end{array} \right.
$$
This form defines a `thick' Reeb foliation $\mathcal{F}$ on $D^2 \times S^1$ (that is, there is a subset $N$ such that $\mathcal{F}|_{N}$ is a Reeb foliation and $\mathcal{F}|_{D^2 \times S^1 \backslash N}$ is diffeomorphic to a product foliation by tori).

Assume that in coordinates $(r, \phi, t)$ the matrix of $g$ has a form:
$$
 g = \left ( \begin{array}{ccc} 1 & 0 & 0 \\
 0 & G(r) & 0 \\
 0 & 0 & 1 \end{array} \right )
$$
In order to calculate the second fundamental form of $\mathcal{F}$ consider the following sections: $X = \frac{\partial}{\partial \phi}, Y =(1 - f(r))\frac{\partial}{\partial r} - f(r) \frac{\partial}{\partial t}$ of the tangent bundle $T\mathcal{F}$. Let $n = f(r) \frac{\partial}{\partial r} + (1 - f(r))\frac{\partial}{\partial t}$ be a normal vector field.

By the straightforward calculation we obtain that the matrix of the second fundamental form is equal to:
$$
\frac{1}{2f(r)^2 - 2f(r) +1}\left ( \begin{array} {cc} - f \frac{\partial G}{\partial r} & 0 \\
                           0 & - (1-f) \frac{\partial f}{\partial r} \end{array} \right )
$$
It is obvious that since $f = 0$ on $[0, \frac{1}{3})$ the foliation by disks is totally geodesic for every choice of $G=G(r)$. Define $G = G(r)$ in the following way:
$$
G = \left \{ \begin{array} {l} r^2, \mbox{when $r \in [0, \frac{1}{4})$} \\
                               \mbox{strictly increasing, when $r \in [\frac{1}{4}, \frac{1}{3})$} \\
                               \mbox{$1$, when $r \in [\frac{1}{3}, 1]$}   \end{array} \ \right.
$$
For this choice of $G$, the metric $g$ is regular in the neighborhood of the core curve $r = 0$ and satisfies conditions of the lemma.
\subsection{Parabolic foliation on $T^2 \times [0,1]$.}
Using the similar arguments we may obtain the following result:
\begin{lem}There exists a foliation $\mathcal{F}$ and metric $g$ on $T^2 \times [0, 1]$ such that
\begin{enumerate}
\item{$\mathcal{F}$ is parabolic with respect to $g$.}
\item{The foliation $\mathcal{F}|_{T^2 \times [0, \frac{1}{3}]}$ is a foliation by the totally geodesic tori $T^2 \times \{r\}$ and the foliation $\mathcal{F}|_{S^1 \times S^1 \times [\frac{2}{3}, 1]}$ is a foliation by the totally geodesic annuli $\{ t \} \times S^1 \times [\frac{2}{3}, 1]$.}
\end{enumerate}
\end{lem}
\emph{Proof:} On $T^2 \times [0,1]$ define the following coordinates:
$$
T^2 \times [0, 1] = \{((\phi, t), r): \ r\in [1, 2], \phi, t \in [0, 2\pi)\}
$$
and consider one-from $\alpha$
$$
 \alpha = f(r)dr + (1 - f(r)) dt
$$
where $f(r)$ is such a smooth function on $[0,1]$ that
$$
 f(r) = \left \{ \begin{array}{l} 1, \ r \in [0, \frac{1}{3}] \\
                                  \mbox{strictly decreasing, when} \   r \in (\frac{1}{3}, \frac{2}{3}] \\
                                  0, \   r \in (\frac{2}{3}, 1]
  \end{array} \right.
$$
This form defines a foliation $\mathcal{F}$ on $T^2 \times [0, 1]$. Similarly to the proof of Lemma $3.3$ we may define a matrix of $g$(
with respect to the coordinates  $((\phi, t), r)$) in the form:
$$
 g = \left ( \begin{array}{ccc} G(r) & 0 & 0 \\
 0 & 1 & 0 \\
 0 & 0 & 0 \end{array} \right )
$$
In order to calculate the second fundamental form of the leaves, consider the following sections of $T\mathcal{F}$: $X = \frac{\partial}{\partial \phi}, Y =(1 - f(r))\frac{\partial}{\partial r} - f(r) \frac{\partial}{\partial t}$. The fiels $n = f(r) \frac{\partial}{\partial r} + (1 - f(r))\frac{\partial}{\partial t}$ is a normal vector field

The second fundamental form of $\mathcal{F}$ with respect to the unit normal $\frac{n}{\abs{n}}$ is given by:
$$
\frac{1}{2f(r)^2 - 2f(r) +1}\left ( \begin{array} {cc} - f \frac{\partial G}{\partial r} & 0 \\
                           0 & - (1-f) \frac{\partial f}{\partial r} \end{array} \right )
$$
It is obvious that since $f = 0$ on $(\frac{2}{3}, 1]$, the foliation $\mathcal{F}_{T^2 \times (\frac{2}{3}, 1]}$ by the horizontal annuli is totally geodesic for an arbitrary choice of $G(r)$. Define $G=G(r)$ by the following formula:
$$
G = \left \{ \begin{array} {l} 1, \mbox{when $r \in [0, \frac{2}{3})$} \\
                               \mbox{strictly decreasing, when $r \in [\frac{2}{3}, \frac{3}{4})$} \\
                               \mbox{strictly increasing, when $r \in [\frac{3}{4}, \frac{4}{5})$} \\
                               \mbox{$r^2$, when $r \in [\frac{4}{5}, 1]$}   \end{array} \ \right.
$$
\begin{figure}
\includegraphics[width = 6.5cm]{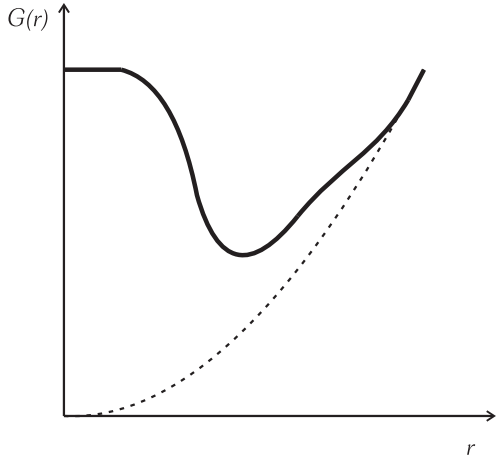}
\caption{The function $G(r)$ in the construction of metric on $T^2 \times [0, 1]$.}
\end{figure}
For this choice of $G$, a metric $g$ satisfies all conditions of the lemma.

\section{Parabolic foliations on three-sphere.}
The aim of this section is to define a parabolic turbularization of $\mathcal{F}$ on $S^3$ along the knot $K$. We will define this foliation in several steps. First, we will define a parabolic Reeb foliation on $S^3$. Then we will construct some special parabolic turbularization of $\mathcal{F}_R$ along the trivial link consisting of $n$ components. For a knot $K$ we will consider its special presentation, which would coincide with the trivial link everywhere except the double points in the frontal projection. In this way we will be able to define a parabolic foliation everywhere except some balls in $S^3$ containing these double points. In order to define a parabolic foliation inside these balls we will `twist' the turbularization along trivial link with two components. Finally, we will show how to glue this foliations back into the sphere to get the desired foliation.

\subsection{Parabolic Reeb foliation on $S^3$.}
\begin{prop}\cite{Bol}
A three-sphere admits a parabolic foliation.
\end{prop}
\emph{Proof:} Take some presentation of a three-sphere $S^3 = D^2_1 \times S^1 \cup_h D^2_2 \times S^1$ as a union of two solid tori. Define the parabolic foliations $\mathcal{F}_1$ and $\mathcal{F}_2$ inside these solid tori as in Lemma $3.3$. In the coordinates $(t, \phi)$ on $\partial D^2 \times S^1$ the gluing diffeomorphism $h$ is given by the matrix $h =\left ( \begin{array}{cc} 0 & 1 \\ -1 & 0 \end{array} \right )$. It is obvious that $h$ is an isometry of the boundary torus $\partial (D^2_1 \times S^1)$. Since the metrics on the solid tori are the direct product metrics and the foliations are the direct product foliations in the (one-sided) neighborhoods of the boundary torus, there is a well-defined glued foliation $\mathcal{F}_R$ and the metric on $S^3$. This foliation is parabolic with respect to a glued metric.
\subsection{Parabolic turbularization along the trivial link.}
\begin{prop}
For every $n \in \mathbb{N}$ there is a parabolic foliation on $S^3$ with $n$ `thick' Reeb components inside the solid torus $D^2(\frac{1}{3}) \times S^1 \subset D^2_1 \times S^1 \subset S^3$.
\end{prop}
\emph{Proof:} Consider a foliation $\mathcal{F}_R$ on $S^3$ defined in proposition $4.1$. First, note that metric inside the disk $D^2(\frac{1}{3}) \times \{0\}$ is a standard euclidian metric. Let $\{x_1, x_2, \ldots, x_n\}$ be a set of vertices of the regular polygon lying inside $D^2(\frac{1}{3}) \times \{ 0 \}$, with the center at $r = 0$ and the radius of the circumscribed circle equal to $\frac{1}{8}$ with respect to a metric induced on the disk $D^2(\frac{1}{3}) \times \{0\}$. Instead of radius $\frac{1}{8}$ we may choose any number such that the circle with the center at the median of the side $\overline{x_ix_{i+1}}$ and the radius equal to its length would entirely lie inside the disk $D^2(\frac{1}{3})$.
Take such $\eps$ that any two circles with the centers at the vertices of the polygon and of radius $\eps$ would be disjoint. Consider the set of vertical circles $\{ x_i \times S^1 \}$ passing through the vertices $x_i$.

On the solid torus $D^2 \times S^1$ take the following coordinates:
$$
D^2 \times S^1 = \{((r, \phi), t): \ r\in [0, \varepsilon], \phi, t \in [0, 2\pi)\}
$$
and consider the function $f$ given by the formula:
$$
f(r) = \left\{ \begin{array}{l} 0, \mbox{when} \ r \in [0, \frac{\eps}{6}] \\
\mbox{strictly increasing, when} \ r \in (\frac{\eps}{6}, \frac{\eps}{3}] \\
1, \mbox{when}\  r \in (\frac{\eps}{3}, \frac{2\eps}{3}] \\
\mbox{strictly decreasing, when} \ r \in (\frac{2\eps}{3}, \frac{5\eps}{6}] \\
0, \mbox{when}\  r \in (\frac{5\eps}{6}, \varepsilon] \\
 \end{array} \right.
$$
The one-form $\alpha = f(r)dr + (1 - f(r))dt$ defines a foliation $\mathcal{F}'$ on $D^2 \times S^1$.

In order to define the metric on $D^2 \times S^1$ consider the following function $G = G(r)$ on it:
$$
G = \left \{ \begin{array} {l} 0, \mbox{when $r \in [0, \frac{\varepsilon}{8})$} \\
                               \mbox{strictly inceasing, when $r \in [\frac{\varepsilon}{8}, \frac{\varepsilon}{6})$} \\
                               \mbox{$\epsilon$,  $r \in [\frac{\varepsilon}{6}, \frac{5\varepsilon}{6})$} \\
                               \mbox{strictly increasing, when$r \in [\frac{5\varepsilon}{6}, \frac{6\varepsilon}{7}]$} \\
                               \mbox{$r^2$, when $r \in [\frac{6\varepsilon}{7}, \varepsilon]$}   \end{array} \ \right.
$$
Define a Riemannian metric $g$ on $D^2 \times S^1$ by the following matrix:
$$
 g = \left ( \begin{array}{ccc} 1 & 0 & 0 \\
 0 & G(r) & 0 \\
 0 & 0 & 1 \end{array} \right )
$$
It is easy to verify that $\mathcal{F}'$ is parabolic with respect to this metric.

Cut $\epsilon$-tubular neighborhoods of the circles $\{x_i\} \times S^1$ and glue solid tori $(D^2 \times S^1, \mathcal{F}')$ instead of each such circle by the identity map. Since $\mathcal{F}'|_{\{T^2 \times [\frac{5\epsilon}{3}, 2\epsilon]\}}$ is a foliation by the totally geodesic annuli it glues correctly to a foliation of $D^2(\frac{1}{3}) \times S^1$ by the horizontal disks.

From the construction of the metric on each $D^2 \times S^1$ it smoothly glues with the (euclidian) metric on $D^2(\frac{1}{3}) \times S^1$ since in the neighborhood of the boundary $\partial(D^2 \times S^1)$ this metric is given by $g = dt^2 + dr^2 + r^2d\phi^2$ (relatively to the coordinate system $((r, \phi), t)$ on $D^2(\frac{1}{3}) \times S^1$). Denote obtained foliation by $\mathcal{F}_n$.
\begin{defn}
We call thus defined foliation $\mathcal{F}_n$ a trivial turbularization with $n$ strings.
\end{defn}
\subsection{Standard presentation of a knot.}
Assume that $K$ is a knot in $S^3$. We may isotope it in such a way that $K$ is a closure of some braid lying inside the solid torus $D^2(\frac{1}{3} )\times S^1$ and transverse to a foliation of $D^2(\frac{1}{3}) \times S^1$ by the totally geodesic disks. Write the presentation of $K$ as a product of transpositions $K = \sigma_1^{\pm 1} \sigma_2^{\pm 1} \ldots \sigma_N^{\pm 1}$.
\begin{figure}
\includegraphics[width = 6.5cm]{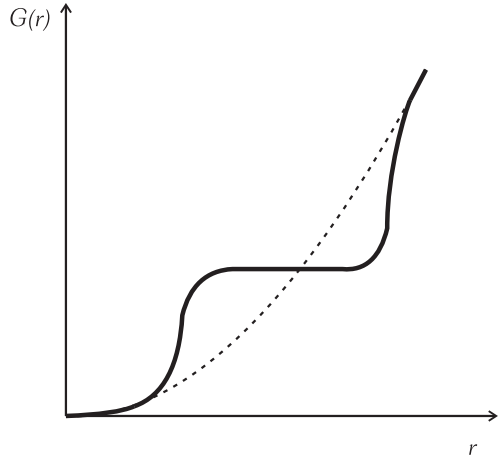}
\caption{Construction of metric on $D^2 \times S^1$.}
\end{figure}
Without loss of generality we may assume that in the frontal projection
$$
 f : D^2 \times S^1 \to [-\frac{1}{3}, \frac{1}{3}] \times S^1 \ \ \ f(x, y, t) = (x, t)
$$
\begin{figure}[htbp]
\includegraphics[width = 5cm]{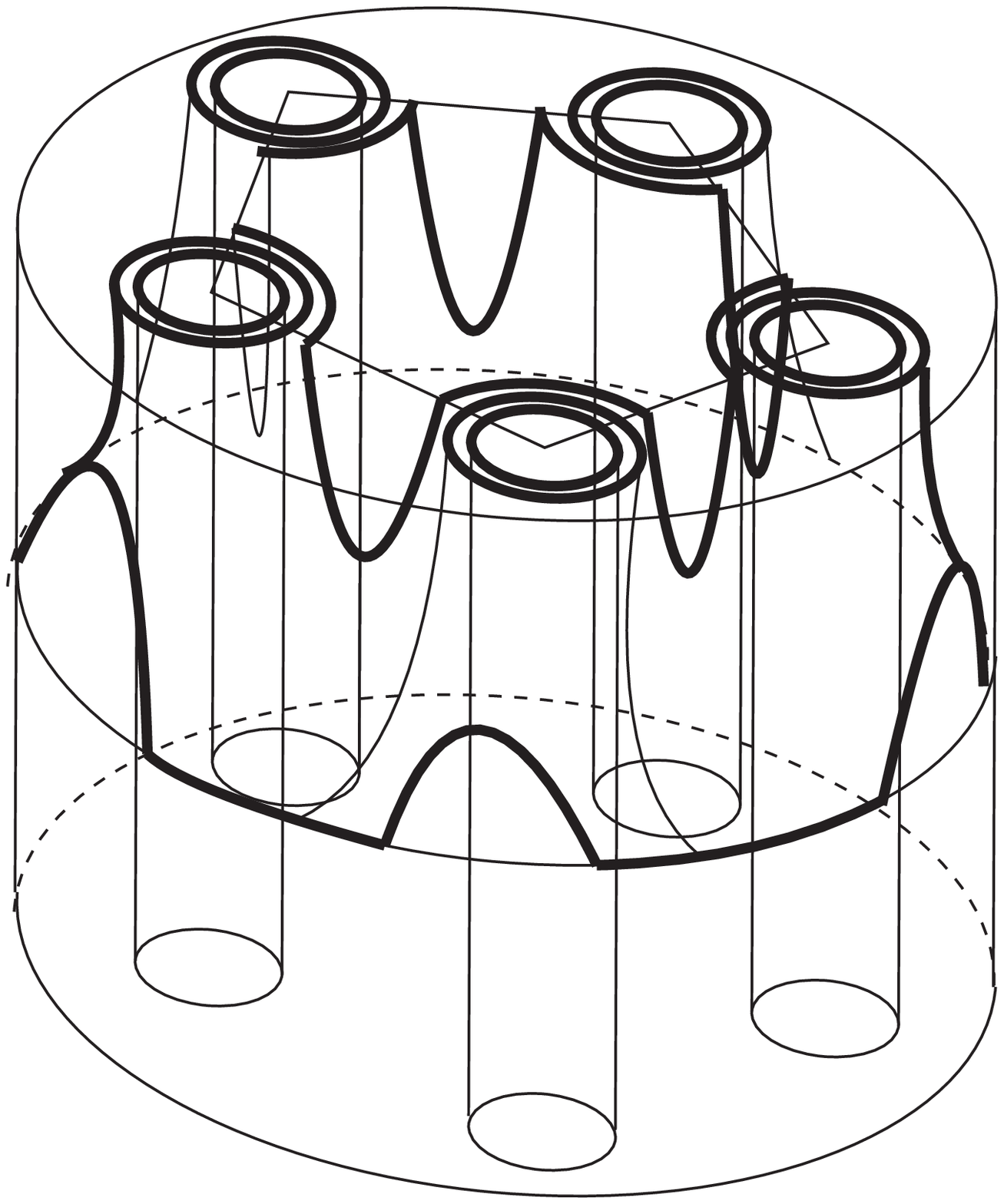} \ \ \ \ \ \ \ \
\includegraphics[width = 5cm]{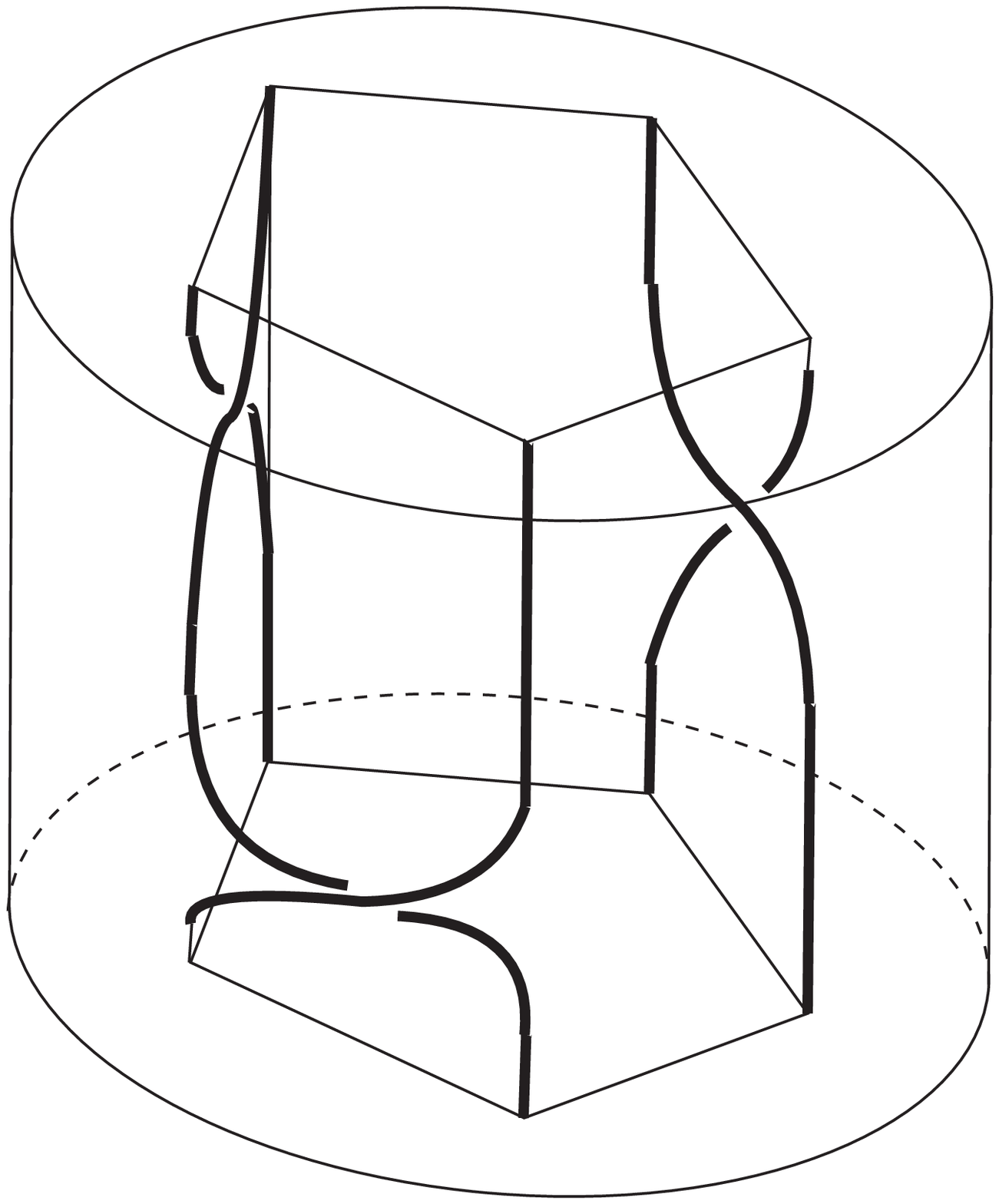}
\caption{Trivial turbularization with 5 strings (left). Standard presentation of a knot(right).}
\end{figure}
there is a finite number of distinct levels $t_1, t_2, \ldots, t_N$ such that $K$ at these points has transverse double points (each point corresponds to a transposition). We can further isotope $K$ in such a way that it would be a subset of $\bigcup_{k = 1}^n \{x_k \times S^1\}$  for some $n$, maybe except the neighborhoods of the inverse images of the double points $f^{-1}(t_i), \ i = 1, 2, \ldots, N$.
\begin{defn}
We call such presentation of $K$ a standard presentation of a knot with $n$ strings.
\end{defn}
Trivial turbularization $\mathcal{F}_n$ with $n$ strings coincides with the turbularization along the standard presentation of $K$ everywhere except the balls around the inverse images of double points.
\subsection{Parabolic foliation in the neighborhood of the transposition.}
In order to define turbularizations along the transpositions consider a trivial turbularization $\mathcal{F}_2$ with two strings on $D^2(\frac{1}{3}) \times [0,1]$ (here we slightly abuse the notation and call trivial turbularization the foliation induced on $D^2(\frac{1}{3}) \times [0,1] = \overline{D^2 \times S^1 \backslash D^2}$). Denote by $g$ the corresponding metric such that $\mathcal{F}_2$ is parabolic. Let $\delta$  be some small real number. Define the following bump function $f = f(r)$ on $[0, \frac{1}{3}]$:
$$
 f(r) = \left \{ \begin{array}{l} \pi, \ r \in [0, \frac{1}{4}] \\
                                  \mbox{strictly decreasing, when} \   r \in (\frac{1}{4}, \frac{1}{3} - \delta] \\
                                  0, \   r \in (\frac{1}{3} - \delta, \frac{1}{3}]
  \end{array} \right.
$$

Consider a one-dimensional dynamical system $(D^2(\frac{1}{3}), \psi_t)$ generated by the flow of the vector field $X = f(r)\frac{\partial}{\partial \phi}$ on the disk $D^2(\frac{1}{3}) = \{(r, \phi): r \in [0,\frac{1}{3}], \phi \in [0, 2\pi)\}$.
\begin{figure}[htbp]
\centering
\includegraphics[width = 6cm]{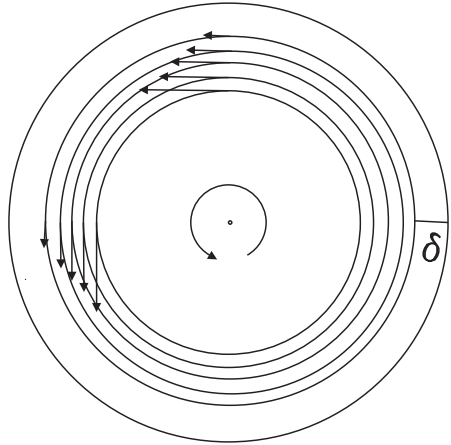}
\caption{Dynamical system in the definition of the left parabolic transposition.}
\end{figure}
Let $h(t)$ be a smooth increasing bump function on $[0, 1]$ such that $h(0) = 0$ and $h(1) = 1$. Associate with it the following diffeomorphism $\Phi$ of $D^2 \times [0, 1]$:
$$
\Phi(r, \phi, t) = (\psi_{h(t)} (r, \phi), t)
$$
\begin{figure}[htbp]
\includegraphics[height =6.5cm, width = 5cm]{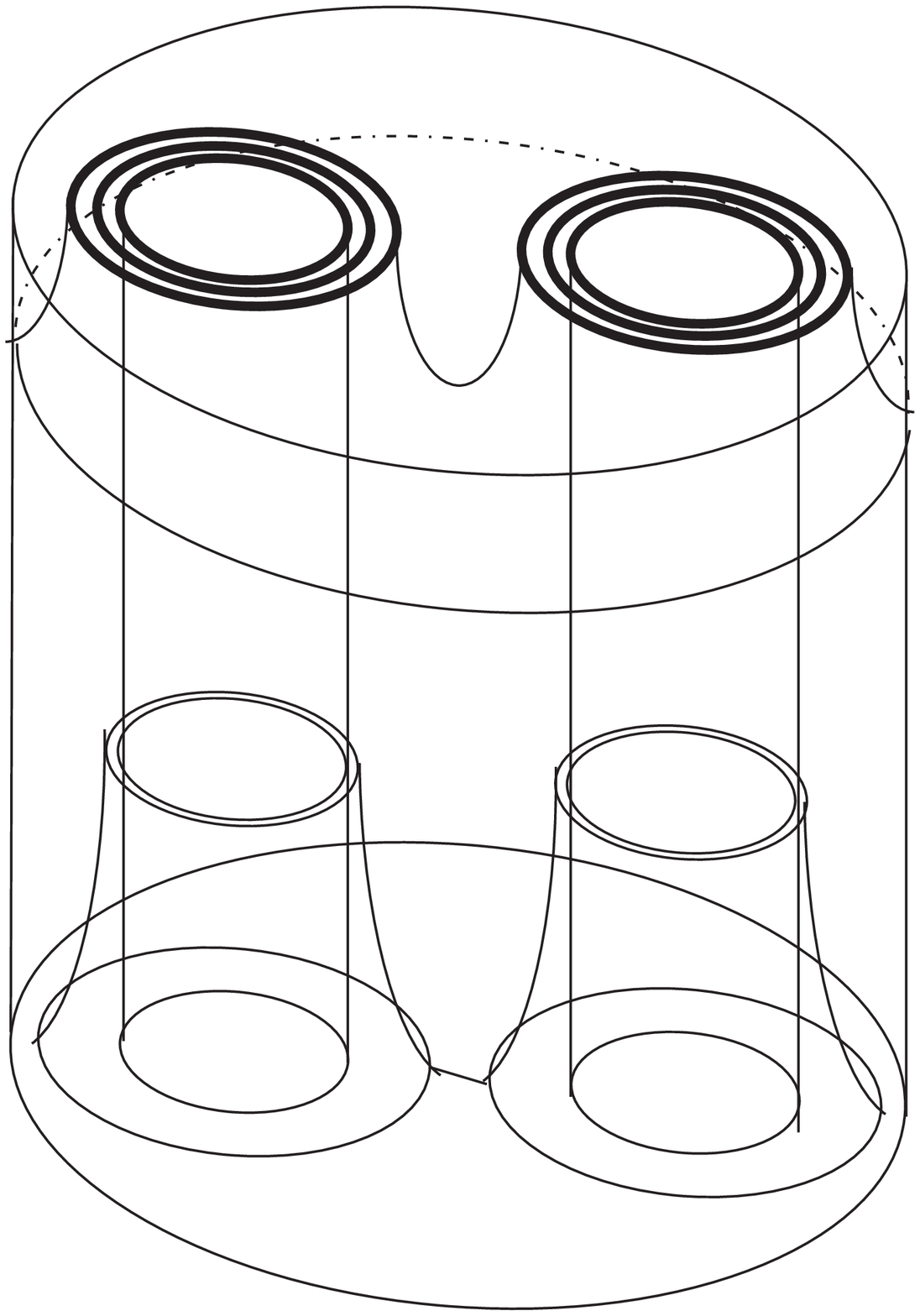} \ \ \ \ \ \ \ \ \ \
\includegraphics[height =6.5cm, width = 5cm]{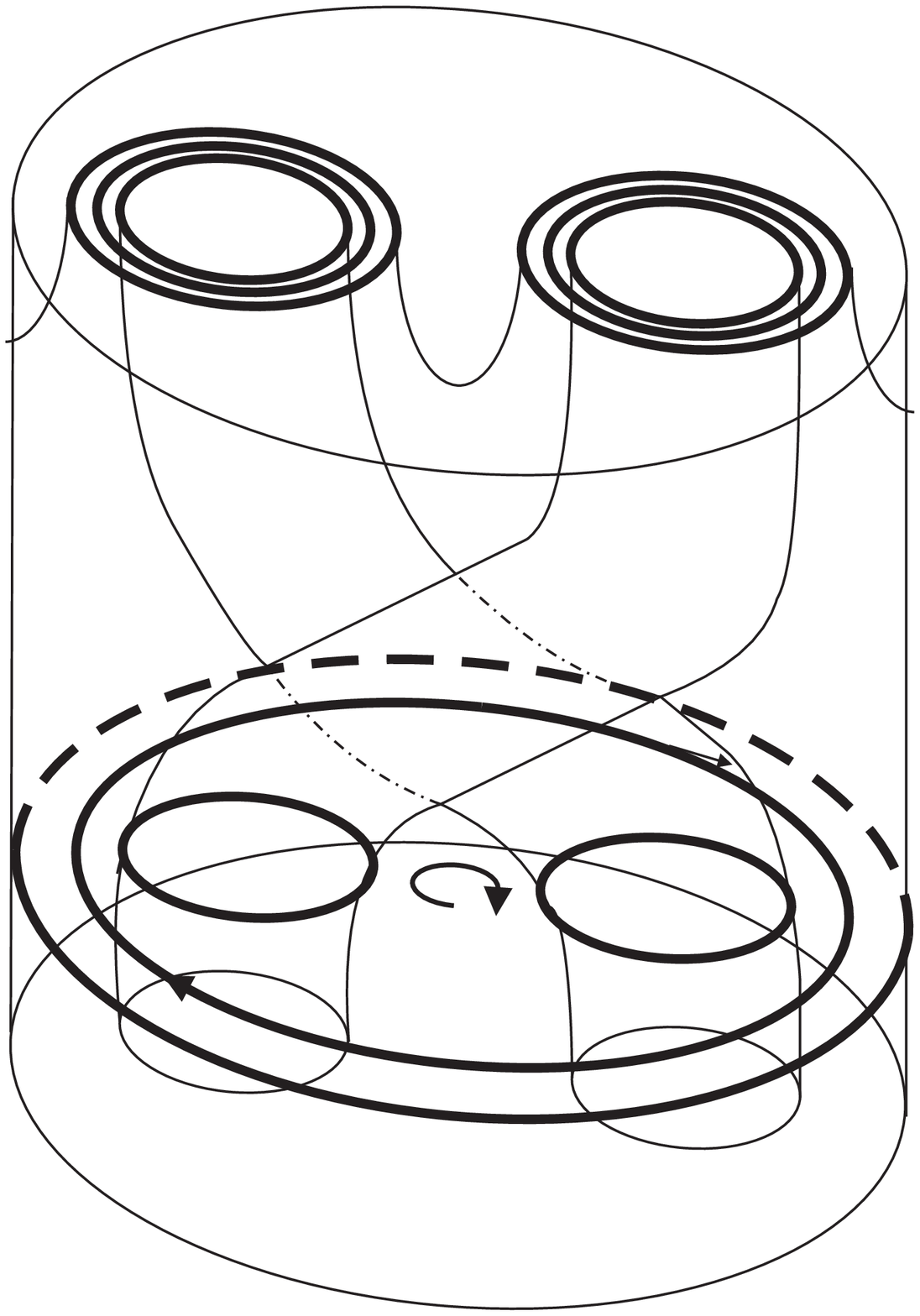}
\caption{Trivial parabolic transposition with two strings(left). Left parabolic transposition(right).}
\end{figure}
This diffeomorphism defines a foliation $\Phi(\mathcal{F}_2)$ on $D^2\times [0, 1]$. This foliation is clearly parabolic with respect to the pull-back metric $(\Phi^{-1})^\ast g$.

Assume that $G(r, \phi)$ is a metric induced on a disc $D^2(\frac{1}{3}) \times \{1\}$ by the metric $g$ and $H(r, \phi)$ -- metric induced on this disk by $(\Phi^{-1})^\ast g$. We may use Lemma $3.1$ to interpolate between these two metrics.

First, note that $G(r, \phi) = H(r, \phi)$ on a disk $D^2(\frac{1}{4}) \times \{1\}$ since $\mathcal{F}_2$ is invariant under the rotation by $\pi$. It is also clear that $G(r, \phi) = H(r, \phi)$ on $(D^2(\frac{1}{3}) \backslash D^2(\frac{1}{3} - \delta)) \times \{1\}$.

Let $N = D^2(\frac{1}{3} - \delta)\backslash D^2(\frac{1}{4}) \times [0, 1]$. From Lemma $3.1$ there is such metric $g$ on $N$ that
\begin{enumerate}
\item{In some tubular neighborhood of $(D^2(\frac{1}{3} - \delta)\backslash D^2(\frac{1}{4})) \times \{0\}$ the metric $g(p,t) = G(p) + dt^2$ for all $p \in D^2(\frac{1}{3} - \delta)\backslash D^2(\frac{1}{4})$}
\item{In some tubular neighborhood of $(D^2(\frac{1}{3} - \delta)\backslash D^2(\frac{1}{4})) \times \{1\}$ the metric $g(p,t) = H(p) + dt^2$ for all $p \in D^2(\frac{1}{3} - \delta)\backslash D^2(\frac{1}{4})$}
\item{$\mathcal{F}$ is parabolic on $N$ with respect to $g$}
\item{There is a neighborhood $U$ of the boundary $\partial (D^2(\frac{1}{3} - \delta)\backslash D^2(\frac{1}{4}))$ such that for all $t \in [0,1], \   g(p, t)|_{\Sigma^2 \times \{t\}} = G(p)$}
\end{enumerate}
On $(D^2(\frac{1}{3}) \backslash D^2(\frac{1}{3} - \delta)) \times [1,2]$ and $D^2(\frac{1}{4}) \times [1, 2]$ consider the direct product foliations.
They are parabolic (even totally geodesic) with respect to a metric $ds^2 = dr^2 + r^2d\phi^2 + dt^2$. Since in the neighborhood of the boundary $\partial N$ the foliation is a direct product foliation and the metric is a direct product metric there is a parabolic foliation correctly defined on the union $L = (D^2(\frac{1}{3}) \backslash D^2(\frac{1}{3} - \delta)) \times [1,2] \cup N \cup D^2(\frac{1}{4}) \times [1, 2]$.

Finally, consider the union $\Phi(\mathcal{F}_2) \cup L$. It is obvious that foliations and metrics on $L$ and $\Phi(\mathcal{F}_2)$ are smoothly glued with each other and define on the union the structure of the parabolic foliation. We are remained to `normalize' this foliation in the $t$ direction. For this consider the map
$$
F : D^2 \times [0, 1] \to \Phi(\mathcal{F}_2) \cup L
$$
defined by the formula $F((r, \phi), t) = ((r, \phi), 2t)$. A foliation formed by the inverse images of the leaves is parabolic in the pull-back metric $F^{\ast}g$.

We call the foliation obtained (together with Riemannian metric) on  $F(\Phi(\mathcal{F}_2) \cup X)$ a standard left (right) parabolic transposition.

\begin{rem}
Note, that we cannot use Lemma $3.1$ directly to the foliation $\Phi(\mathcal{F}_2)$ since this foliation is not a foliation by disks.
\end{rem}
\subsection{Parabolic turbularization along $K$ on a three-sphere $S^3$.}
\begin{lem}
For each topological knot type $K$ there is a parabolic foliation $\mathcal{F}_K$ on $S^3$ such that $\mathcal{F}_K$ is a parabolic turbularization along $K$ of the parabolic Reeb foliation on $S^3$ (see proposition $4.1$).
\end{lem}
\emph{Proof:} Consider a standard presentation of $K$. Let $n$ be a number of strings in it. Recall, that $K$ is a subset of the union of vertical circles $\{x_i \times S^1\}$ everywhere except some neighborhoods of the inverse images of the double points of $K$ in the frontal projection.

In the frontal projection, let $t_1, \ldots, t_N$ denote the set of $t$-coordinates of the double points of $K$.

\begin{figure}
\centering
\includegraphics[height =6.5cm, width = 6.5cm]{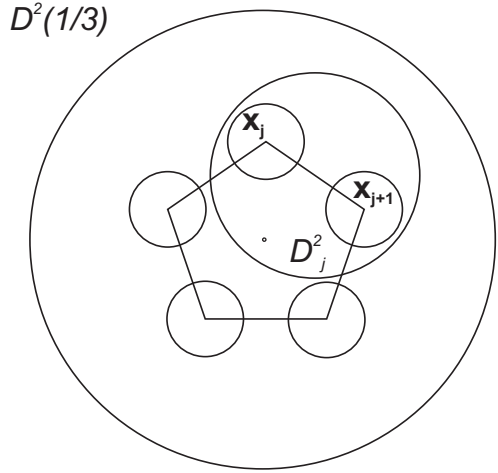}
\caption{A disk $D^2(\frac{1}{3}) \times \{0\}.$}
\end{figure}

Write the presentation of $K$ as a product $K = \sigma_1^{\pm 1} \sigma_2^{\pm 1} \ldots \sigma_N^{\pm 1}$. Recall that with each vertex of the regular polygon we associated the disk with the center at the vertex and with radius $\eps$(see page $15$). For each $\sigma_j$ consider a disk $D^2_j$ with the center at the median of the edge $\overline{x_{j-1}x_j}$ and with radius $d(\overline{x_{j-1}x_j})/ 2 + 2\eps$. We choose such disk to make sure that the small disks with the centers at the vertices and of radius $\eps$ are lying inside it(see Figure. 9).

Note, that since $x_j$ are the vertices of regular polygon the points $x_{j-1}$ and $x_j$ are the only points from the set $\{x_1, x_2, \ldots, x_n\}$, which lie inside the disk $D^2_j$. Consider a set of disjoint intervals $I(t_j), j = 1\ldots N \subset [0, 1]$ such that $t_j \in I(t_j)$.

On $S^3$ define a standard turbularization with $n$ strings and a metric (see Proposition $4.2$) such that $\mathcal{F}_n$ is parabolic with respect to it. For each $j$ let  $(D^2(\frac{1}{3}) \times [0, 1], \mathcal{F}', g)$ be a left (or right) standard parabolic transposition depending on the degree of the corresponding $\sigma_j$. Denote the radius of $D^2_j$ by $r_j$ and the length of the segment $I(t_j)$ by $d_j$. Consider the map $$F_j: D^2(\frac{1}{3}) \times [0, 1] \to D^2_j \times I(t_j)$$ which is given by the formula
$$
F_j((r, \phi), t) = ((\frac{3r}{r_j}, \phi), \frac{t}{d_j})
$$
This map defines a foliation $F_j(\mathcal{F}')$ inside the ball $D^2_j \times I(t_j)$. An inverse map $F^{-1}_j$ defines on $D^2_j \times I(t_j)$ such metric that $F_j(\mathcal{F}')$ is parabolic with respect to it. Since in the neighborhood of the gluing the foliation is a direct product foliation and the metric is a direct product metric it glues correctly to $\mathcal{F}_n$. Denote obtained foliation by $\mathcal{F}_K$. This foliation is parabolic with respect to the glued metric.

\section{Gluing the solid torus.}
\emph{Proof of theorem $1.1$:}
In order to finish the proof of theorem we have to perform a Dehn surgery on a knot $K$.

Consider the foliation $\mathcal{F}_K$ on $S^3$. Let $N$ denote such tubular neighborhood of $K$ that $\partial N = T^2$ is a leaf of $\mathcal{F}_K$. Denote by $X = \overline{S^3 \backslash N}$ and consider an arbitrary diffeomorphism $f$:
$$
f : \partial X \to \partial(D^2 \times S^1)
$$
Up to isotopy this diffeomorphism is defined by the map it induces in homology:
$$
 f_{\ast} : H_1(T^2) \to H_1(T^2) \ \ \ f_{\ast} \in SL_2(\mathbb{Z})
$$
In particular we may think that $f = \left ( \begin{array}{cc} a & b \\ c & d \end{array} \right )$ is a linear map.

On $D^2 \times S^1$ define a parabolic foliation and a metric as in Lemma $3.3$. This metric is euclidian in some neighborhood of the boundary torus. Therefore $f$ defines the following metric on $\partial X$:
$$
 G  = \left ( \begin{array}{cc} a^2 + c^2 & ac + bd \\  ac+ bd & b^2 + d^2  \end{array} \right )
$$
In its own part on $\partial X$ we have a metric $H = \left (\begin{array}{cc}1 & 0 \\ 0 & 1\end{array}\right )$. In order to interpolate between $G$ and $H$ consider the union $X\cup T^2\times [0, 1]\cup_f D^2 \times S^1$. From the Corollary $3.2$ there is a metric on $T^2 \times [0,1]$ which deforms $G$ to $H$ in such a way that the foliation by tori $T^2 \times \{pt\}$ is parabolic. Since in (one-sided) neighborhoods of $\partial X$ and $\partial (D^2 \times S^1)$ the metrics are the direct product metrics, on the union  $X\cup T^2\times [0, 1]\cup_f D^2 \times S^1$ we obtain smooth Riemannian metric. Therefore on $X\cup T^2\times [0, 1]\cup_f D^2 \times S^1 \sim X\cup_f D^2 \times S^1 $ and therefore on every closed oriantable three-manifold we were able to define a parabolic foliation.

\end{document}